\documentclass[a4paper,12pt,reqno]{amsart}

\usepackage{amsmath,amssymb,amsthm,amsaddr}
\usepackage[margin=3.4cm,top=3.5cm,bottom=3.8cm,footskip=1cm,headsep=1.5cm]{geometry}
\usepackage[utf8]{inputenc}
\usepackage{graphicx}
\usepackage{xcolor}
\usepackage{doi}

\usepackage[noend]{algpseudocode}
\usepackage[ruled]{algorithm2e}
\usepackage{booktabs}
\usepackage{siunitx}
\usepackage{bm}
\usepackage[list=true, font=small, labelfont=bf, 
labelformat=brace, position=top]{subcaption}

\usepackage{enumitem}

\def\div{\operatorname{div}}
\def\curl{\operatorname{curl}}
\def\hh{\bm{H}}

\def\bb{\bm{B}}

\def\jj{\bm{J}}

\def\ee{\bm{E}}

\theoremstyle{definition}

\usepackage{url}

\begin{document}
\title[Magnetic Space-Time Formulation Including Hysteresis]{Inclusion of an Inverse Magnetic Hysteresis Model into the Space-Time Finite Element Method for Magnetoquasistatics}

\author{M.~Gobrial$^1$, L.~Domenig$^2$, M.~Reichelt$^1$, \\
M.~Kaltenbacher$^2$, O.~Steinbach$^1$}

\address{%
\small 
$^1$Institute of Applied Mathematics, TU Graz, Austria \\
$^2$Institute of Fundamentals and Theory in
Electrical Engineering, TU Graz, Austria}

\begin{abstract}
In this note we discuss the numerical solution of the
eddy current approximation of the Maxwell equations using
the simple Pragmatic Algebraic Model to include hysteresis
effects. In addition to the more standard time-stepping
approach we propose a space-time finite element method
which allows both for parallelization and adaptivity
simultaneously in space and time. Numerical experiments confirm
both approaches yield the same numerical results.
\end{abstract}

\maketitle 

\begin{quote}
\end{quote}

\bigskip

\section{Introduction} \label{sec:intro}
For the mathematical modeling of the electromagnetic behavior 
in an electric machine, e.g. an electric motor or a transformer, 
we consider the eddy current approximation of Maxwell's equations in the low-frequency regime, cf., e.g.,
\cite{AlonsoValli2010, Jackson1999, Monk2003},
\begin{equation}\label{eq:Maxwell}
\curl \hh = \mathbf{j}_s, \quad
\div \bb = 0, \quad
\curl \ee = -  \partial_t \bb \quad
\mbox{in} \; D \subset {\mathbb{R}}^3,
\end{equation}
with boundary conditions
$\bb \cdot \mathbf{n} = 0$ on $\partial D$.
For simplicity we consider a fixed domain $D$,
for which we can write Ohm's law as
$\jj = \sigma \ee$; for the case of a moving domain, see, e.g.,
\cite{Gangl:2025}. In addition to the Maxwell equations
\eqref{eq:Maxwell} we need to have a constitutive relation
$\hh = \hh(\bb)$ which in many cases is described in terms of 
the magnetic reluctivity $\nu$ and, in presence of permanent magnets, the permanent magnetization $\mathbf{M}$, i.e.
\begin{align}
\label{eq:BH_nu}
    \hh = \mathbf{\nu}(\lVert \bb \rVert) \, \bb - \mathbf{M}.
\end{align}
Note that $\nu$ may depend on the magnitude of the magnetic 
flux density $\bb$ as e.g. for ferromagnetic materials, see, e.g., 
\cite{Meunier2008}. Although the constitutive law 
\eqref{eq:BH_nu} covers a wide range of physically relevant phenomena, it neglects the effects of hysteresis. These effects
become more and more important to accurately describe electrical
devices, such as electric machines and transformers. For this reason, it is necessary to employ material models that cover hysteresis into the mathematical model, and its numerical
simulation by finite element methods. For the standard finite element method, this has already successfully been done in the static case using the magnetic scalar potential \cite{Domenig2024_incorphyst} as well as using the magnetic vector potential \cite{xiao2022modeling}. For transient simulations, there exist also approaches to consider the hysteresis effect in the finite element method using a magnetic vector potential \cite{Jacques2018,purnode2024neural}. To the authors knowledge, however, there exists no approach to incorporate hysteretic material models into a space-time finite element framework.

Since the magnetic flux density is a solenoidal vector field, 
$\bb = \curl \mathbf{A}$, we can rewrite the eddy current problem
\eqref{eq:Maxwell} in its equivalent vector potential formulation \cite{Monk2003, Kaltenbacher2015}
\begin{equation}\label{eq:EC_vector_potential}
    \sigma \partial_t \mathbf{A} + \curl [\nu \curl \mathbf{A} ]
    = \mathbf{j}_s + \curl \mathbf{M},
\end{equation}
where we have used the constitutive law \eqref{eq:BH_nu}. Usually, gauging techniques are applied in order to ensure the uniqueness of the vector potential $\mathbf{A}$, see, e.g.,
\cite{Kaltenbacher2015, Hafner1987}. A common simplification 
of the eddy current problem \eqref{eq:EC_vector_potential} is the reduction to a two-dimensional model problem assuming that one dimension of the computational domain is much larger than the other dimensions and that the geometry is invariant in the larger scale, cf.~\cite{Clain1993}. Hence, we can pose the eddy current problem on a cross section $\Omega \subset {\mathbb{R}}^2$ of the 
electric machine in $D$, where the electromagnetic quantities 
take the form
\begin{align*}
    \hh = \begin{pmatrix}
        H_1(x_1,x_2,t) \\ H_2(x_1,x_2,t) \\ 0
    \end{pmatrix}, \
    \mathbf{M} = \begin{pmatrix}
        M_1(x_1,x_2,t) \\ M_2(x_1,x_2,t) \\ 0
    \end{pmatrix}, \
    \mathbf{j}_s = \begin{pmatrix}
        0 \\ 0 \\ j_s(x_1,x_2,t)
    \end{pmatrix}.
\end{align*}
It follows that $\mathbf{j}_s$ is divergence free by construction, and the magnetic flux density $\bb$ must admit the same form as the electric field intensity $\hh$, due to~\eqref{eq:BH_nu}. 
Using $\bb = \curl \mathbf{A}$ we further have
\begin{equation}
\label{eq:A=u}
    \mathbf{A} = \begin{pmatrix}
        0 \\ 0 \\ u(x_1,x_2,t)
    \end{pmatrix}, \quad \bb = \begin{pmatrix}
        \partial_{x_2} u(x_1,x_2,t) \\ -\partial_{x_1} u(x_1,x_2,t) \\ 0
    \end{pmatrix}.
\end{equation}
With this we can rewrite \eqref{eq:EC_vector_potential} 
for $(x,t) \in Q := \Omega \times (0,T)$ as
\begin{equation}
\label{eq:EC_2D}
    \sigma(x) \partial_t u(x,t) - \div_x \big[
    \nu(x,|u|) \nabla_x u(x,t)\big] = j_s(x,t) - 
    \div_x M^\perp(x,t),
\end{equation}
where $M^\perp = (-M_2,M_1)^\top$ is the perpendicular of the magnetization and $T>0$ is the final time. For the sake of completeness, boundary conditions $u=0$ on 
$\Sigma := \partial \Omega \times (0,T)$
as well as an initial condition $u(0)=0$ in $\Omega$
must be set. Although, this approach~\eqref{eq:EC_2D} of the eddy current problem is indeed a typical approximation in the context of electric machines, it does not consider hysteresis effects. In this work, we will apply a specific law among many, that considers the effect of hysteresis in terms of the constitutive law~\eqref{eq:BH_nu}.

The remainder of this work is organized as follows. In Section~\ref{sec:hysteresis_model} we introduce the Pragmatic Algebraic Model (PAM), a specific vector hysteresis model previously analyzed 
in~\cite{xiao2022modeling, purnode2024neural, Henrotte2015}. 
We incorporate this constitutive law into Maxwell's equations 
\eqref{eq:Maxwell} and formulate the problem on the 
cross section $\Omega$ of the geometry $D$ to derive the underlying parabolic evolution problem \eqref{eq:EC_hyst}. 
Section~\ref{sec:fem} presents two numerical approaches for solving the resulting time-dependent partial differential equation. On the one hand, we formulate a classical time stepping method~\cite{Thomee2006}, and on the other hand, we describe a space-time finite element method motivated by~\cite{Steinbach2015}. A comparative analysis of both is provided in Section~\ref{sec:numerical_results}, indicating a good agreement of the results.

\section{Hysteresis Model}\label{sec:hysteresis_model}
The aim of this section is to modify the eddy current problem
\eqref{eq:EC_2D} to a model, which can consider hysteresis effects. A physical system subject to hysteresis does not only depend on the input data but also on the history of these
data \cite{Brokate1996}. For ferromagnetic materials, e.g., iron, hysteresis effects are quite natural and shall be considered in the constitutive law.

In this work we consider the Pragmatic Algebraic Model (PAM) as  hysteresis model, which uses six real positive parameters 
$p_j \in {\mathbb{R}}_+$, $j=0,\ldots,5$,
in one algebraic expression in order to describe hysteresis, see~\cite{Xiao2020}. In contrast to other hysteresis models, the efficiency of PAM lies in its formulation by one algebraic expression, which takes the static and the dynamic effects into account. The adapted constitutive law reads as
\begin{equation}\label{eq:hysteresis_law}
    \hh = f(\bb) \bb + g(\partial_t{\bb})\partial_t{\bb} - \mathbf{M},
\end{equation}
where
\[
    f(\bb) = p_0 + p_1\lVert\bb\rVert^{2p_2}, \quad
    g(\partial_t{\bb}) = p_3 + \frac{p_4}{\sqrt{p_5^2 + \lVert \partial_t{\bb} \rVert^2}}.
\]
The first expression $f(\bb)$ of \eqref{eq:hysteresis_law} describes the anhysteretic part, which is in fact similar to the magnetic reluctivity $\nu$, which reflects the BH-curve relation described by \eqref{eq:BH_nu}. The parameters $p_0, p_1, p_2$ can be fitted in order to obtain the same behavior as $\nu$ in the classical approach. The second expression $g(\partial_t{\bb})$ describes on one hand the macroscopic eddy currents by the parameter $p_3$, and on the other hand the hysteresis effects, that are considered by $p_4$ and $p_5$, cf.~\cite{Xiao2020}. As before, $\mathbf{M}$ is the permanent magnetization of occurring permanent magnets. Now, when using the constitutive law~\eqref{eq:hysteresis_law} instead of~\eqref{eq:BH_nu} and again the vector potential 
$\bb = \curl \mathbf{A}$, the underlying eddy current equation considering hysteresis reads as
\begin{equation}
\label{eq:EC_vector_potential_hyst}
    \sigma \partial_t \mathbf{A} + \curl\bigg( f\big(\curl(\mathbf{A})\big) \curl(\mathbf{A}) + g({\curl(\partial_t\mathbf{A})}\big){\curl(\partial_t \mathbf{A})} \bigg)
    = \mathbf{j}_s + \curl(\mathbf{M}).
\end{equation}
The reduction to the spatially two-dimensional case requires the same 
assumptions as above, hence the vector potential~$\mathbf{A}$ has the same form as 
in~\eqref{eq:A=u}, and we can rewrite~\eqref{eq:EC_vector_potential_hyst}
as
\begin{equation}\label{eq:EC_hyst}
    \sigma \, \partial_t u - \div_x 
    \big[ f(|\nabla_x u|) \nabla_x u + g(|\partial_t \nabla_x u|) \partial_t \nabla_x u \big] 
    = j_s - \div_x M^\perp 
\end{equation}
in $Q := \Omega \times (0,T)$.
In addition to the partial differential equation \eqref{eq:EC_hyst} we consider homogeneous Dirichlet boundary 
conditions $u=0$ on $\Sigma := \partial\Omega \times (0,T)$, which implies the induction boundary 
condition $B \cdot n = 0$, and the initial 
condition $u(x,0) = 0$ for $x \in \Omega$.

\section{Finite Element Formulation} \label{sec:fem}

\subsection{Time Stepping Framework}
When multiplying the time dependent partial differential equation
\eqref{eq:EC_hyst} with a spatial test function $v$ vanishing
on $\partial \Omega$, integrating over $\Omega$ and
applying integration by parts, we obtain
\begin{eqnarray}\label{eq:semi VF}
\int_\Omega \sigma \, \partial_t u \, v \, dx
+ \int_\Omega \Big[ f(|\nabla_x u|) \, \nabla_x u 
+ g(|\partial_t \nabla_x u|) \,
\partial_t \nabla_x u \Big] \cdot \nabla_x v \, dx
\hspace*{2cm} && \\ \nonumber = \int_\Omega \Big[ j_s \, v 
+ M^\perp \cdot \nabla_x v \Big] \, dx \, .
\end{eqnarray}
Let $S_h^1(\Omega) = \mbox{span} \{ \phi_k \}_{k=1}^{M_\Omega}$
be the standard finite element space of piecewise linear basis
functions $\phi_k$ which are defined with respect to an
admissible decomposition of the computational domain 
$\Omega$ into shape regular triangular finite elements
$\tau_\ell$ of the spatial mesh size $h_x$, $\ell=1,\ldots,N_\Omega$, and which are zero on $\partial \Omega$.
The semi-discretization
of \eqref{eq:semi VF} is then equivalent of a system of
nonlinear ordinary differential equations,
\begin{equation}\label{eq:semi NLGS}
[ M_h + A_h(\dot{u}_h)]\underline{\dot{u}}(t) 
+ K_h(u_h) \underline{u}(t) =
\underline{F}(t),
\end{equation}
where the entries of the mass and stiffness matrices 
as well as of the load vector are
given by, for $j,k=1,\ldots,M_\Omega$,
\begin{eqnarray*}
M_h[j,k] & = & \int_\Omega \sigma(x) \, \phi_k(x) \,
\phi_j(x) \, dx, \\
K_h(u_h)[j,k] & = & \int_\Omega f(|\nabla_x u_h(x,t)|) 
\, \nabla_x \phi_k(x) \cdot \nabla_x \phi_j(x) \, dx, \\
A_h(\dot{u}_h)[j,k] & = & \int_\Omega g(|\nabla_x \dot{u}_h(x,t)|) \,
\partial_t \nabla_x \phi_k(x) \cdot \nabla_x \phi_j(x) \, dx, \\
F_j(t) & = & \int_\Omega \Big[ j_s(x,t) \, \phi_j(x) 
+ M^\perp(x,t) \cdot \nabla_x \phi_j(x) \Big] \, dx.
\end{eqnarray*}
Note that $\underline{u}(t) \in {\mathbb{R}}^{M_\Omega}$
is the vector of the time dependent coefficients of the
numerical solution
\[
u_h(x,t) = \sum\limits_{k=1}^{M_\Omega} u_k(t) \phi_k(x) \, .
\]
For time discretization we introduce a temporal mesh size $h_t$
and we define time steps $t_i = i h_i$, $i=0,\ldots,N_T$.
When considering
\[
u_h(x,t_i) = \sum\limits_{k=1}^{M_\Omega} u_k(t_i) \phi_k(x)
= \sum\limits_{k=1}^{M_\Omega} u_k^i \phi_k(x) =
u_h^i(x),
\]
and using the backward finite difference scheme
\[
\dot{u}_h(x,t_i) \simeq 
\frac{1}{h_t} [u_h(x,t_i)-u_h(x,t_{i-1})] =
\frac{1}{h_t} \sum\limits_{k=1}^{M_\Omega}
[u_k^i-u_k^{i-1}] \phi_k(x),
\]
the time discretization of \eqref{eq:semi NLGS} results in a
sequence of nonlinear systems of algebraic equations,
$i=1,\ldots,N_T$,
\begin{equation}\label{eq:NLGS time}
\frac{1}{h_t} [ M_h + A_h([u_h^i-u_h^{i-1}]/h_t)]
[\underline{u}^i - \underline{u}^{i-1}] 
+ K_h(u_h^i) \underline{u}^i(t) =
\underline{F}(t_i),
\end{equation}
with the initial condition $\underline{u}^0 = \underline{0}$.
At each time step $t_i$, $i=1,\ldots,N_T$, the nonlinear system
\eqref{eq:NLGS time} is solved via Newton's method where we
use the results as given in \cite{purnode2024neural} to
compute all involved derivatives analytically.

\subsection{Space-Time Framework}
Next we consider a space-time variational formulation for the eddy current problem\eqref{eq:EC_hyst}. We now multiply the
transient partial differential equation \eqref{eq:EC_hyst}
with a test function $v(x,t)$ vanishing on 
$\Sigma = \partial \Omega \times (0,T)$, and integrate
over the space-time domain $Q = \Omega \times (0,T)$. Integration by 
parts only with respect to the spatial components finally provides the space-time variational formulation
\begin{align}
\label{eq:VF_space_time_direct}
    & \int_0^T \int_\Omega \sigma \partial_t u \, v \, dx \, dt
    + \int_0^T \int_\Omega f(|\nabla_x u|) \nabla_x u \cdot \nabla_x v \, dx \, dt \\
    & \hspace*{5mm}  
    + \int_0^T \int_\Omega g(|\partial_t \nabla_x u|) \partial_t \nabla_x u \cdot \nabla_x v \, dx \, dt 
    = \int_0^T  \int_\Omega \Big[ j_s \, v 
    + M^\perp \cdot \nabla_x v \Big] dx \, dt . 
    \nonumber
\end{align}
Let $S_h^1(Q) = \mbox{span} \{ \varphi_k \}_{k=1}^{M_Q}$ be the
space-time finite element space of piecewise linear basis
functions $\varphi_k$ which are defined with respect to an
admissible decomposition of the space-time domain $Q$
into tetrahedral finite elements $q_\ell$, $\ell=1,\ldots,N_Q$
of mesh size $h$, and which are zero at initial time $t=0$,
and on the lateral boundary $\Sigma$.
However, since second order derivatives occur in the weak 
formulation \eqref{eq:VF_space_time_direct}, we can not use
$S_h^1(Q)$ for a conforming finite element discretization
of \eqref{eq:VF_space_time_direct}. Instead, we use the
substitution $p(x,t) := \partial_t u(x,t)$ to rewrite
\eqref{eq:VF_space_time_direct} as
\begin{align}
\label{eq:VF_space_time_system_1}
    & \int_0^T \int_\Omega \sigma \partial_t u \, v \, dx \, dt
    + \int_0^T \int_\Omega f(|\nabla_x u|) \nabla_x u \cdot \nabla_x v \, dx \, dt \\
    & \hspace*{5mm}  
    + \int_0^T \int_\Omega g(|\nabla_x p|) \nabla_x p \cdot \nabla_x v \, dx \, dt 
    = \int_0^T  \int_\Omega \Big[ j_s \, v 
    + M^\perp \cdot \nabla_x v \Big] dx \, dt ,
    \nonumber
\end{align}
together with a second variational formulation
\begin{align}
\label{eq:VF_system}
    \int_0^T \int_\Omega p(x,t) q(x,t) \, dx \, dt 
    = \int_0^T \int_\Omega \partial_t u(x,t) q(x,t) \, dx \, dt .
\end{align}
Since $p = \partial_tu$ only has zero boundary conditions on
$\Sigma$, there are no initial conditions at $t=0$.
For the space-time finite element approximation of $p$ we
therefore have to used the extended finite element space
$\widetilde{S}_h^1(Q) = 
\mbox{span} \{ \varphi_k \}_{k=1}^{\widetilde{M}_Q}$,
where the additional basis functions $\varphi_k$,
$k=M_Q+1,\ldots,\widetilde{M}_Q$ are related to the nodes
at $t=0$. The space-time finite element discretization
of the system \eqref{eq:VF_space_time_system_1} and
\eqref{eq:VF_system} is then equivalent to a system of
nonlinear algebraic equations,
\begin{equation}\label{eq:space-time NLGS}
    \left( \begin{array}{cc}
B_h + K_h(u_h) & A_h(p_h) \\[1mm]
- \widetilde{B}_h & M_h
    \end{array} \right)
    \left( \begin{array}{c}
    \underline{u} \\[1mm]
    \underline{p} 
    \end{array} \right)
    =
    \left( \begin{array}{c}
    \underline{F} \\[1mm] \underline{0}
    \end{array} \right),
\end{equation}
where the entries of the block matrices are given by,
for $k,\ell=1,\ldots,M_Q$ and $i,j=1,\ldots,\widetilde{M}_Q$,
\begin{eqnarray*}
B_h[\ell,k] 
& = & \int_0^T \int_\Omega \sigma(x) \, \partial_t \varphi_k(x,t)
\, \varphi_\ell(x,t) \, dx \, dt, \\
K_h(u_h)[\ell,k]
& = & \int_0^T \int_\Omega f(|\nabla_x u_h(x,t)|) \,
\nabla_x \varphi_k(x,t) \cdot \varphi_\ell(x,t) \, dx \, dt, \\
A_h(p_h)[\ell,i]
& = & \int_0^T \int_\Omega g(|\nabla_x p_h|) \,
\nabla_x \varphi_i(x,t) \cdot \nabla_x \varphi_\ell(x,t)
\, dx \, dt, \\
M_h[j,i] & = & \int_0^T \int_\Omega \varphi_i(x,t) \cdot
\varphi_j(x,t) \, dx \, dt , \\
\widetilde{B}_h[j,k]
& = & \int_0^T \int_\Omega \partial_t \varphi_k(x,t) \,
\varphi_j(x,t) \, dx \, dt \, .
\end{eqnarray*}
In addition, the entries of the load vector are given as
\[
F_\ell = \int_0^T \int_\Omega \Big[ j_s(x,t) \, \varphi_\ell(x,t)
+ M^\perp(x,t) \cdot \nabla_x \varphi_\ell(x,t) \Big] \, dx \, dt,
\quad \ell=1,\ldots,M_Q,
\]
and $\underline{u} \in {\mathbb{R}}^{M_Q}$ and
$\underline{p} \in {\mathbb{R}}^{\widetilde{M}_Q}$ are the
coefficient vectors
of the finite element functions $u_h$ and
$p_h$, respectively. Since the space-time mass matrix $M_h$
is invertible, we can compute 
$\underline{p} = M_h^{-1} \widetilde{B}_h \underline{u}
\leftrightarrow p_h$ to conclude the nonlinear 
Schur complement system
\begin{equation}\label{eq:nonlinear Schur}
\Big[ B_h + K_h(u_h) + A_h(p_h) M_h^{-1} \widetilde{B}_h
\Big] \underline{u} = \underline{F}.
\end{equation}
To solve the global nonlinear system \eqref{eq:nonlinear Schur}, 
we use an exact Newton method with Armijo's damping 
strategy~\cite{Deuflhard2011}, where the linearized system of each Newton iteration is solved with the parallel direct solver MUMPS supported by PETCs~\cite{Dalcin2011}, which is based on a mesh decomposition method provided by the finite element library Netgen/NGSolve~\cite{Schoberl1997}. 

\section{Numerical Results}
\label{sec:numerical_results}
In this section we want to compare the proposed methods. While the space-time formulation~\eqref{eq:nonlinear Schur} has to deal with 
a three-dimensional problem, which is further blown up to a system for the additional variable $p_h$, the time-stepping method simply 
considers the two-dimensional spatial problem. However, the large system of the space-time method needs to be solved in parallel only once, whereas the time-stepping method needs to solve the spatial problem sequentially for each time step. The upcoming examples will 
first give a comparison between these two methods with respect to an academic model problem. Secondly, the Team problem 32~\cite{teamproblem32} will be considered, which verifies the applicability of the proposed hysteresis 
problem~\eqref{eq:hysteresis_law} as well as the accuracy of the introduced methods for solving the eddy current equation including hysteresis.

\subsection{Simple Geometry}
\begin{figure}
\centering
\begin{subfigure}{.5\textwidth}
  \centering
  \includegraphics[width=\textwidth, trim=4cm 0cm 4cm 1cm, clip]{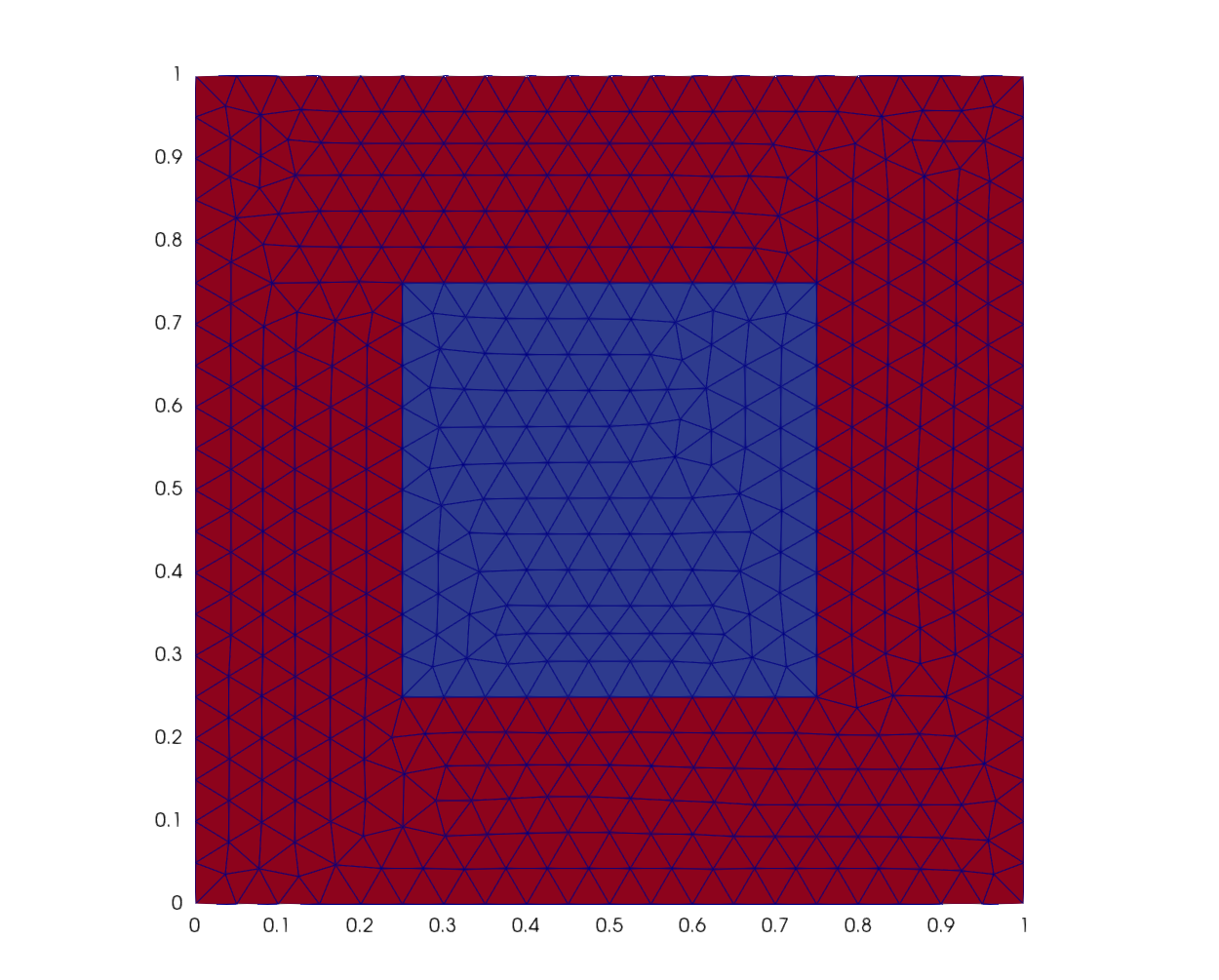}
\end{subfigure}%
\begin{subfigure}{.5\textwidth}
  \centering
  \includegraphics[width=\textwidth, trim=3cm 0cm 0cm 5cm, clip]{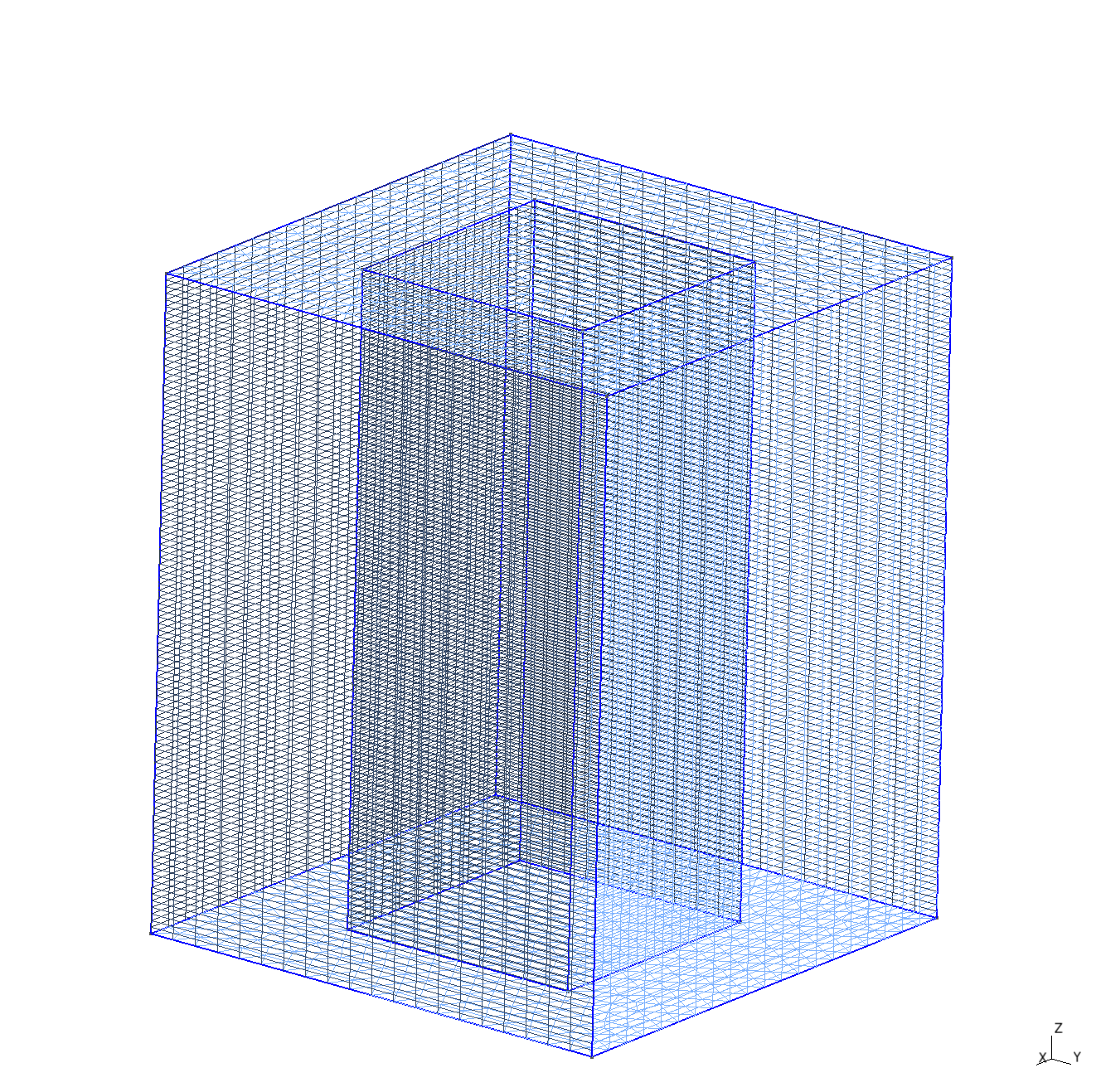}
\end{subfigure}
\caption{Left: The domain $\Omega = (0,1)^2$ consisting of two materials~$\Omega_{cu}$ (in blue) and $\Omega_{fe}$ (in red). Right: The space-time mesh $Q = \Omega \times (0,T)$, which has 100 time slices in temporal direction, 53.530 nodes and 293.400 elements.}
\label{fig:design_cuboid}
\end{figure}

\begin{figure}
\centering
\begin{subfigure}{.5\textwidth}
  \centering
  \includegraphics[width=\textwidth]{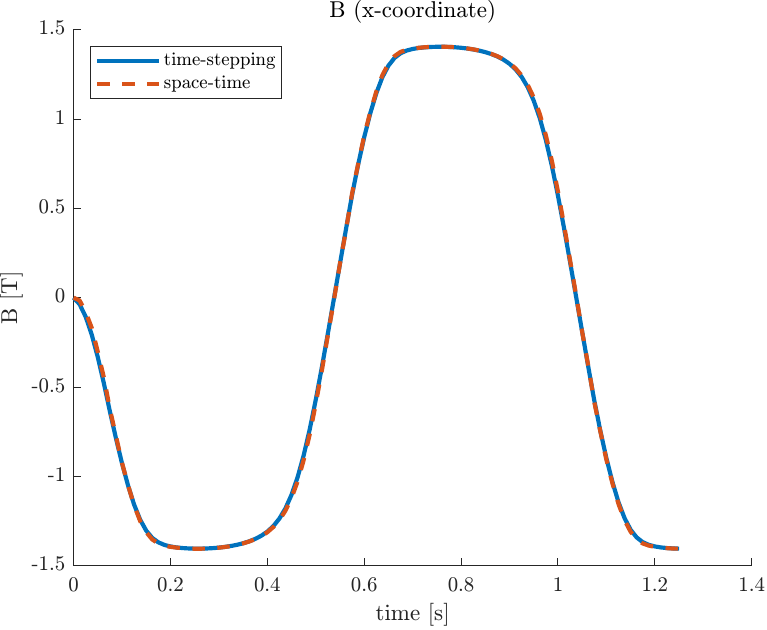}
\end{subfigure}%
\begin{subfigure}{.5\textwidth}
  \centering
  \includegraphics[width=\textwidth]{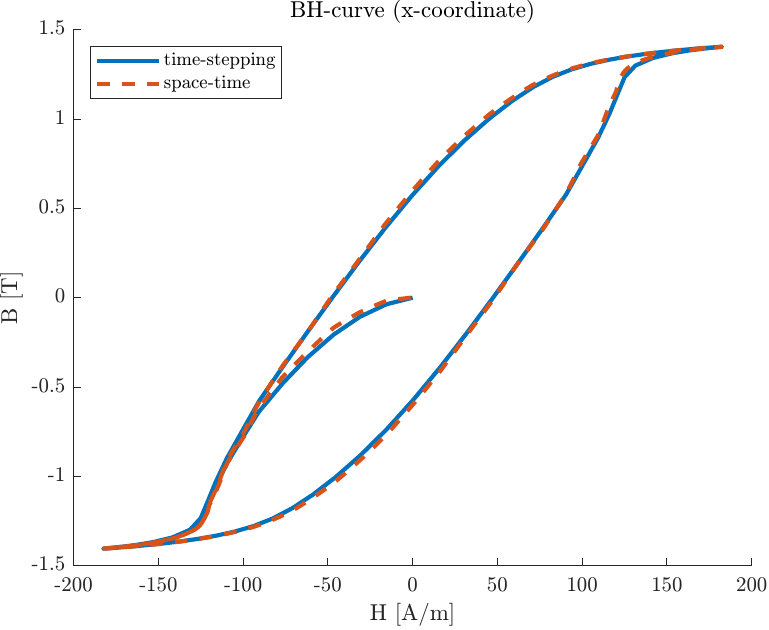}
\end{subfigure}
\caption{Left: The magnetic flux density $B_x$ over time. Right: The BH-curve indicating the hysteresis effect.}
\label{fig:comp_BH_toy_prob}
\end{figure}

The first example considers a two-dimensional spatial domain 
$\Omega = (0,1)^2$, which consists of two different materials, 
$\Omega_{cu} = (0,25,0.75)^2$ consisting of copper through which the excitation $j_s$ passes, and $\Omega_{fe} = \Omega \setminus \overline{\Omega_{cu}}$ consisting of iron in which the hysteresis model is obtained. The final time is given as $T=1.25$. Figure~\ref{fig:design_cuboid} shows the spatial domain, as well as the structured space-time mesh for the space-time finite element method. Furthermore, we use the following parameters,
\begin{align*}
    \sigma(x) &= 
    \begin{cases}
        0 & \text{in } \Omega_{cu}, \\
        0.01 & \text{in } \Omega_{fe},
    \end{cases} 
    & &f(|\nabla_x u|) = 
    \begin{cases}
       \frac{10^7}{4\pi} & \text{in } \Omega_{cu}, \\
       p_0 + p_1 |\nabla_x u|^{2p_2} & \text{in } \Omega_{fe},
    \end{cases} \\
    j_s(x,t) &= \begin{cases}
        2000\sin(2\pi t)  & \text{in } \Omega_{cu}, \\
        0 & \text{in } \Omega_{fe},
    \end{cases}
    & &g(|\partial_t \nabla_x u|) = \begin{cases}
        0 & \text{in } \Omega_{cu}, \\
        p3 + \frac{p_4}{\sqrt{p_5^2 + |\partial_t \nabla_x u|^2}} & \text{in } \Omega_{fe},
    \end{cases}
\end{align*}
where $p_0 = 75.6$, $p_1= 0.0223$, $p_2 = 11.47$, $p_3= 0.0001$, $p_4= 65.8$, $p_5=1$. Note that, the equivalence 
$g(|\nabla_x p|) = g(|\partial_t \nabla_x u|)$ holds on the continuous level, and that $M^\perp = (0,0)^\top$, since no permanent magnets occur.
Figure~\ref{fig:comp_BH_toy_prob} visualizes the magnetic flux density~$\bb_x$ (in the x-component) as well as the hysteresis curve and yields that both methods produce almost the
same results.

\subsection{Two-Phase Transformer -- TEAM Problem 32}
\begin{figure}
\centering
\begin{subfigure}{.4\textwidth}
  \centering
  \includegraphics[width=\textwidth, trim=3cm 3cm 4cm 3cm, clip]{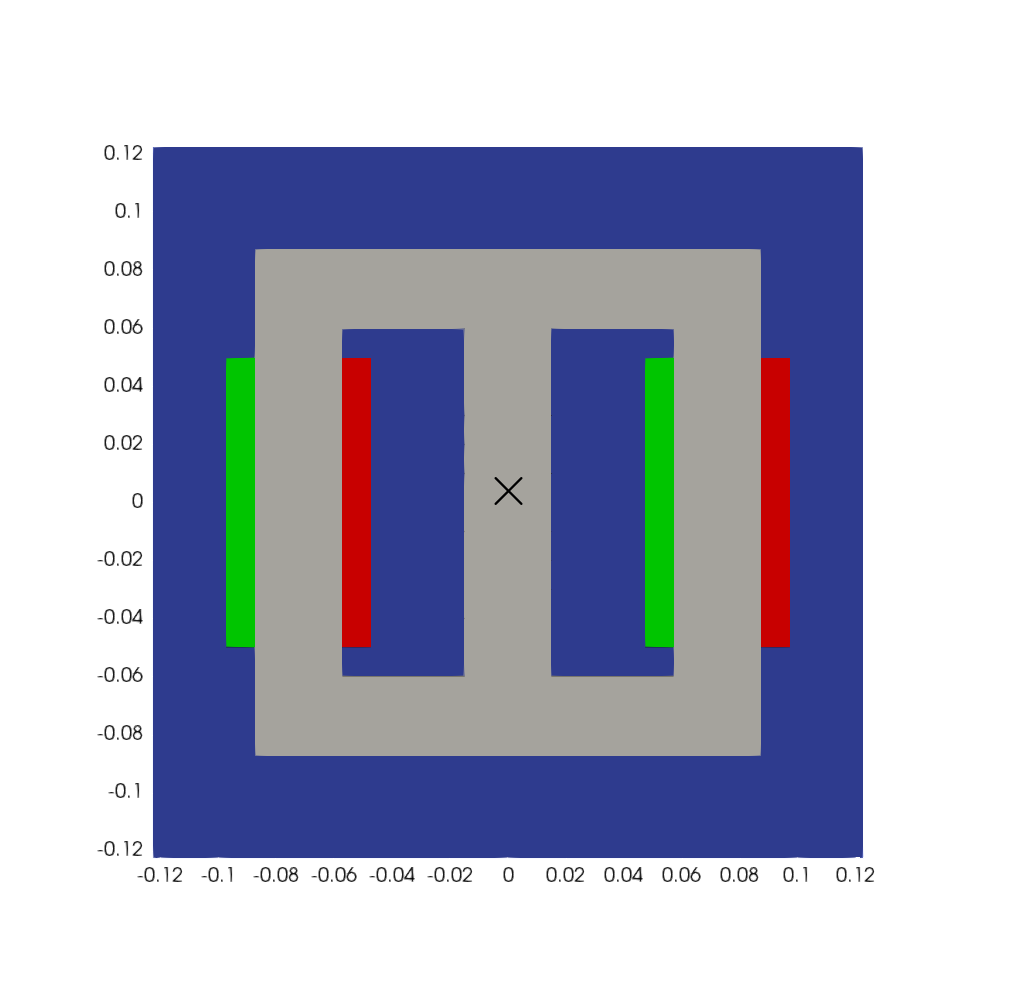}
\end{subfigure}%
\begin{subfigure}{.6\textwidth}
  \centering
  \includegraphics[width=\textwidth, trim=0cm 0cm 0cm 0.25cm, clip]{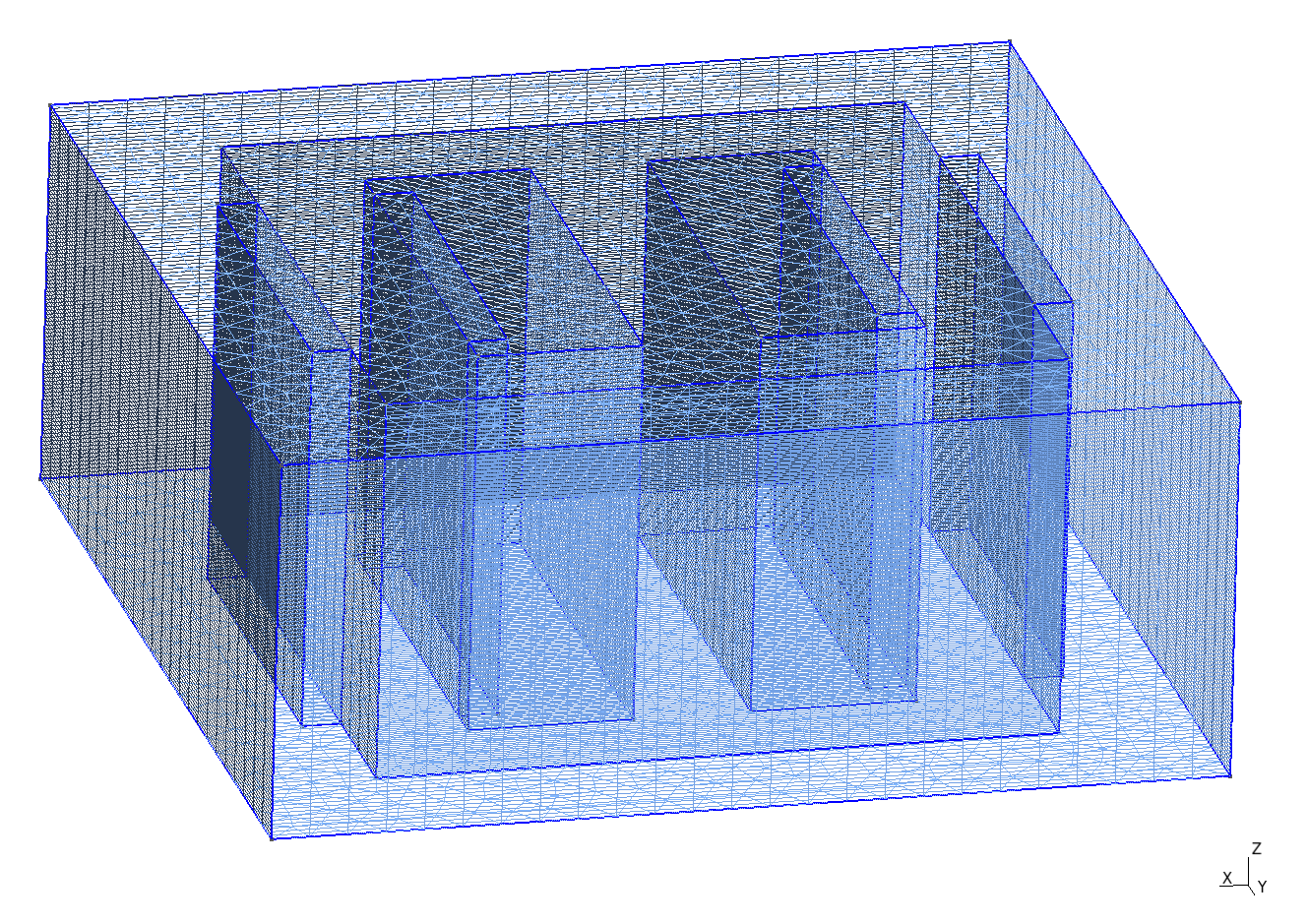}
\end{subfigure}
\caption{Left: The domain $\Omega$ of the transformer consisting of the iron core $\Omega_{fe}$ (in gray), the windings $\Omega_{cu}$ (in green and red) and air $\Omega_a$ (in blue). Right: The space-time mesh $Q = \Omega \times (0,T)$, which has 100 time slices in temporal direction, 91.405 nodes and 206856 elements.}
\label{fig:design_trafo}
\end{figure}

\begin{figure}
\centering
\begin{subfigure}{.5\textwidth}
  \centering
  \includegraphics[width=\textwidth]{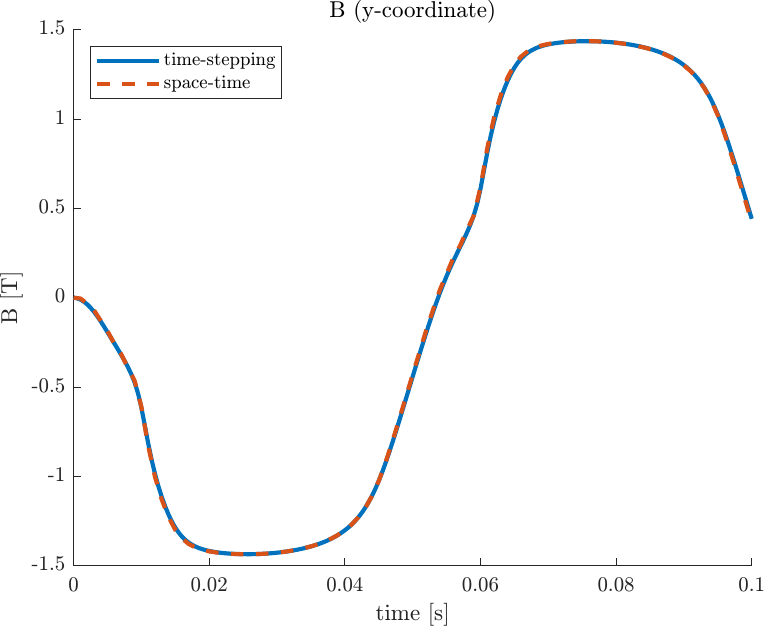}
\end{subfigure}%
\begin{subfigure}{.5\textwidth}
  \centering
  \includegraphics[width=\textwidth]{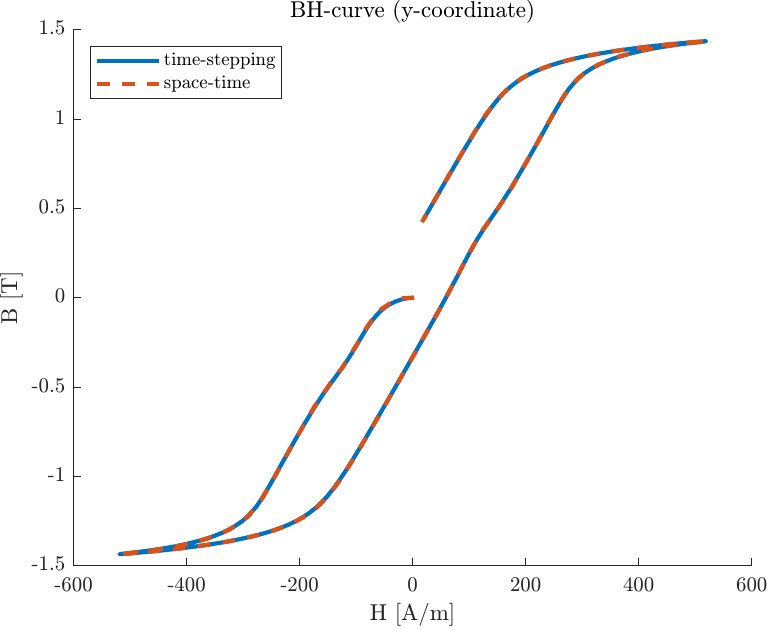}
\end{subfigure}
\caption{Left: The magnetic flux density $B_x$ over time. Right: The BH-curve indicating the hysteresis effect.}
\label{fig:comp_BH_teamProblem}
\end{figure}

\begin{figure}
    \centering
    \includegraphics[width=\textwidth]{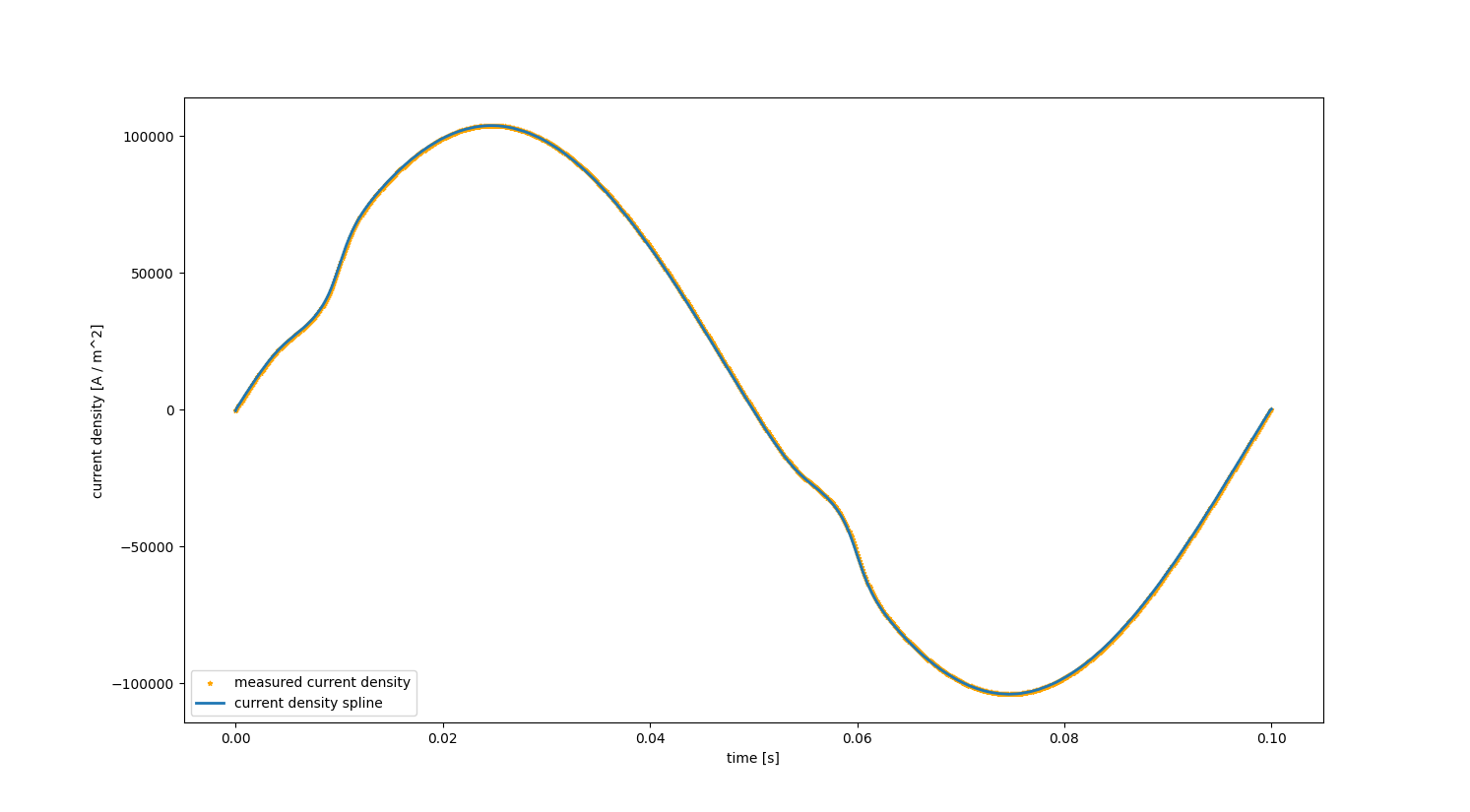}
    \caption{The current density spline $\widetilde{j}$ (in blue) interpolated with respect to the measured values (in orange) of test case 1A from~\cite{teamproblem32}.}
    \label{fig:current_density_trafo}
\end{figure}

Our second example is a two-phase transformer given by the TEAM problem 32~\cite{teamproblem32}, which presents the electric field 
simulation including hysteresis on a three limbed ferromagnetic core with two thin windings, see Figure~\ref{fig:design_trafo}. The two-dimensional computational domain $\Omega$ has the same dimensions as in~\cite{teamproblem32} and consists of three different materials, the iron core $\Omega_{fe}$, the windings $\Omega_{cu}$ at the external limbs and air~$\Omega_{a}$. The considered time span is for $T=0.1$ and the parameters are
\begin{align*}
    \sigma(x) &= 
    \begin{cases}
        0 & \text{in } \Omega_{cu} \cup \Omega_a, \\
        0 & \text{in } \Omega_{fe},
    \end{cases} & & f(|\nabla_x u|) = 
    \begin{cases}
       \frac{10^7}{4\pi} & \text{in } \Omega_{cu} \cup \Omega_a \\
       p_0 + p_1 |\nabla_x u|^{2p_2} & \text{in } \Omega_{fe},
    \end{cases} \\
    j_s(x,t) &= \begin{cases}
        \widetilde{j}(t)  & \text{in } \Omega_{cu}, \\
        0 & \text{in } \Omega_{fe} \cup \Omega_a,
    \end{cases} & &  
    g(|\partial_t \nabla_x u|) = \begin{cases}
        0 & \text{in } \Omega_{cu} \cup \Omega_a \\
        p3 + \frac{p_4}{\sqrt{p_5^2 + |\partial_t \nabla_x u|^2}} & \text{in } \Omega_{fe},
    \end{cases}
\end{align*}
where $p_0 = 181,88232$, $p_1= 0.267053$, $p_2 = 8.999565$, $p_3= 0.00001$, $p_4=  0.0001$, $p_5=50$ and the current density $\widetilde{j}$ is a B-spline interpolation of the measured current values from~\cite{teamproblem32} multiplied with $90$, the number of turns, and divided by the area of the winding, cf. Figure~\ref{fig:current_density_trafo}. Since no permanent magnets occur, we have $M^\perp = (0,0)^\top$. Figure~\ref{fig:comp_BH_teamProblem} depicts the magnetic flux density~$\bb_y$ (in the y-component), which shows, that both methods agree well. It also visualizes the BH-curve, in which the hysteresis effect is visible and both methods indicate this behavior very well.

\section{Conclusion}
In this paper, we have formulated the eddy current approximation derived from Maxwell’s equations in the low-frequency regime, incorporating a specific hysteresis model for the numerical simulation of electromagnetic fields. The nonlinear hysteresis model PAM is represented by a single algebraic expression, making it an efficient choice for capturing hysteresis effects.
To solve the resulting nonlinear time-dependent PDE, we employed two different numerical approaches. The first is a classical semi-discretization method, where the finite element method is used for spatial discretization, followed by an implicit time-stepping scheme. The second approach is a space-time finite element method, which requires solving a saddle-point system but allows the entire problem to be solved at once, enabling parallel computations in spatial and temoral directions simultaneously.
Finally, we compared both methods and demonstrated the applicability of the hysteresis model by achieving similar simulation results. However, the main advantage of space-time
finite element methods is the possibility to use adaptivity
simultaneously in space and time in order to resolve the
potential $u$ locally, and therefore reduce the total
number of degrees of freedom to reach a prescribed accuracy.
A more detailed numerical analysis of such an adaptive
space-time finite element approach will be done in future research.

\begingroup
\footnotesize
\subsection*{Acknowledgements}
This work was supported by the joint DFG/FWF Collaborative Research Centre CREATOR (DFG: Project-ID 492661287/TRR 361; FWF: 10.55776/F90) at TU Darmstadt, TU Graz and JKU Linz.
\endgroup

\bibliographystyle{ieeetr}
\bibliography{refs}

\begin{thebibliography}{10}

\bibitem{AlonsoValli2010}
A.~Alonso~Rodríguez and A.~Valli, {\em Eddy Current Approximation of Maxwell
  Equations: Theory, Algorithms and Applications}.
\newblock Milano: Springer, 2010.

\bibitem{Jackson1999}
J.~Jackson, {\em Classical Electrodynamics}.
\newblock Wiley, 1999.

\bibitem{Monk2003}
P.~Monk, {\em Finite Element Methods for Maxwell's Equations}.
\newblock Oxford Academic, 2003.

\bibitem{Gangl:2025}
P.~Gangl, M.~Gobrial, and O.~Steinbach, ``A space-time finite element method
  for the eddy current approximation of rotating electric machines,'' {\em
  Comput. Methods Appl. Math.}, vol.~25, pp.~441--457, 2025.

\bibitem{Meunier2008}
G.~Meunier, {\em The Finite Element Method for Electromagnetic Modeling}.
\newblock Wiley, 2008.

\bibitem{Domenig2024_incorphyst}
L.~D. Domenig, K.~Roppert, and M.~Kaltenbacher, ``Incorporation of a 3-d
  energy-based vector hysteresis model into the finite element method using a
  reduced scalar potential formulation,'' {\em IEEE Trans. Magn.}, vol.~60,
  no.~6, pp.~1--8, 2024.

\bibitem{xiao2022modeling}
X.~Xiao, F.~M{\"u}ller, M.~M. Nell, and K.~Hameyer, ``Modeling anisotropic
  magnetic hysteresis properties with vector stop model by using finite element
  method,'' {\em COMPEL}, vol.~41, no.~2, pp.~752--763, 2022.

\bibitem{Jacques2018}
K.~Jacques, {\em Energy-Based Magnetic Hysteresis Models - Theoretical
  Development and Finite Element Formulations}.
\newblock PhD thesis, Université de Liège, 2018.

\bibitem{purnode2024neural}
F.~Purnode, F.~Henrotte, G.~Louppe, and C.~Geuzaine, ``Neural network-based
  simulation of fields and losses in electrical machines with ferromagnetic
  laminated cores,'' {\em Int. J. Numer. Model.}, vol.~37, no.~2, p.~e3226,
  2024.

\bibitem{Kaltenbacher2015}
M.~Kaltenbacher, {\em Numerical Simulation of Mechatronic Sensors and
  Actuators}.
\newblock Berlin, Heidelberg: Springer, 2015.

\bibitem{Hafner1987}
C.~Hafner, {\em Numerische {B}erechnung elektromagnetischer {F}elder}.
\newblock Berlin: Springer, 1987.

\bibitem{Clain1993}
S.~Clain, J.~Rappaz, M.~Swierkosz, and R.~Touzani, ``Numerical modeling of
  induction heating for two-dimensional geometries,'' {\em Math. Models Methods
  Appl. Sci.}, vol.~3, no.~6, pp.~805--822, 1993.

\bibitem{Henrotte2015}
F.~Henrotte, S.~Steentjes, K.~Hameyer, and C.~Geuzaine, ``{Pragmatic two-step
  homogenisation technique for ferromagnetic laminated cores},'' {\em IET Sci.
  Meas. Technol.}, vol.~9, no.~2, pp.~152--159, 2015.

\bibitem{Thomee2006}
V.~Thom\'ee, {\em Galerkin finite element methods for parabolic problems},
  vol.~25 of {\em Springer Series in Computational Mathematics}.
\newblock Berlin: Springer, 2006.

\bibitem{Steinbach2015}
O.~Steinbach, ``Space-time finite element methods for parabolic problems,''
  {\em Comput. Methods Appl. Math.}, vol.~15, no.~4, pp.~551--566, 2015.

\bibitem{Brokate1996}
M.~Brokate and J.~Sprekels, {\em Hysteresis and Phase Transitions}.
\newblock Applied Mathematical Sciences, New York: Springer, 1996.

\bibitem{Xiao2020}
X.~Xiao, F.~Müller, G.~Bavendiek, and K.~Hameyer, ``Analysis of vector
  hysteresis models in comparison to anhysteretic magnetization model,'' {\em
  Eur. Phys. J. Appl. Phys.}, vol.~91, no.~2, p.~20901, 2020.

\bibitem{Deuflhard2011}
P.~Deuflhard, {\em Newton Methods for Nonlinear Problems: Affine Invariance and
  Adaptive Algorithms}.
\newblock Springer Series in Computational Mathematics, Springer, 2011.

\bibitem{Dalcin2011}
L.~D. Dalcin, R.~R. Paz, P.~A. Kler, and A.~Cosimo, ``Parallel distributed
  computing using {P}ython,'' {\em Adv. Water Resour.}, vol.~34, no.~9,
  pp.~1124--1139, 2011.

\bibitem{Schoberl1997}
J.~Schöberl, ``{NETGEN}: An advancing front 2d/3d-mesh generator based on
  abstract rules,'' {\em Comput. Visual. Sci.}, vol.~1, pp.~41--52, 1997.

\bibitem{teamproblem32}
O.~Bottauscio, M.~Chiampi, C.~Ragusa, L.~Rege, and M.~Repetto, ``A test case
  for validation of magnetic field analysis with vector hysteresis,'' {\em IEEE
  Trans. Magn.}, vol.~38, pp.~893--896, 2002.

\end{thebibliography}


\newpage 

\appendix

\end{document}